\documentclass[12pt]{article}
\usepackage{amsmath, amsfonts, amssymb, euscript}
\usepackage{graphics}
\input{amssym.def}
\def\Bbb{\mathbb}
\def\goth{\mathfrak}
\newtheorem{theorem}{Theorem}

\date{}

\begin{document}

\begin{center}
\Large{\bf Berezin quantization and representation theory}
\end{center}

\begin{center}
\copyright \ \ {\bf{V.~F.~Molchanov}}
\end{center}

\sloppy

{\footnotesize\begin{center}
Tambov State University named after G.~R.~Derzhavin\\
Internatsionalnaya 33, Tambov, Russian Federation, 392000 \\
E-mail: v.molchanov@mail.ru
\end{center}}

\null

\hfill {\it to memory of Gerrit van Dijk}

\null

{\hspace{0.22cm} {\parbox{14.0cm} {\small \noindent
\noindent \textbf{Abstract.} {We present an approach to Berezin quantization (a variant of quantization in the spirit of Berezin) on para-Hermitian symmetric spaces using the notion of an "overgroup". This approach gives  covariant and contravariant symbols and the Berezin transform in a natural and transparent way}

}}}

\null

\null

The basic notion of the Berezin quantization on a manifold $M$ is a correspondence which to an operator $A$ from a class assigns the pair of functions $F$ and $F^\natural$ defined on $M$. These functions are called {\it covariant and contravariant symbols} of $A$. We are interested in homogeneous space $M=G/H$ (namely, a para-Hermitian symmetric space) and classes of operators related to the representation theory. 

The most algebraic version of quantization -- we call it the {\it polynomial quantization} -- is obtained when operators belong to the algebra of operators corresponding in a representation $T$ of $G$ to elements $X$ of the universal enveloping algebra ${\rm Env} \, (\goth g)$ of the Lie algebra $\goth g$ of $G$. In this case symbols turn out to be polynomials on $\goth g$.

Now we offer a new theme in the Berezin quantization on $G/H$: as an initial class of operators we take operators corresponding to elements {\it of the group $G$} itself in a representation $T$ of this group. We call it the {\it group quantization}.

\null

\setcounter{theorem}{0}
\setcounter{section}{1}
\setcounter{equation}{0}
\setcounter{lemma}{0}

\null

\begin{center}
{{\large {\bf \S \ 1. Para-Hermitian symmetric spaces}}}
\end{center}

\null

Let $G/H$ be a {\it semisimple symmetric space}. Here $G$ is a connected semisimple Lie group with an involutive automorphism $\sigma \ne 1$, and $H$ is an open subgroup of $G^{\sigma}$, the subgroup of fixed points of $\sigma$. We consider that groups act on their homogeneous spaces {\it from the right}, so that $G/H$ consists of right cosets $Hg$.

Let $\goth g$ and $\goth h$ be the Lie algebras of $G$ and of $H$ respectively. Let $B_{\goth g}$ be the Killing form of $G$. There is a decomposition of $\goth  g$ into direct sums of $+1,\,-1$-eigenspaces of the involution
$\sigma$:
$$
{\goth  g} = {\goth  h} + {\goth  q}.
$$
The subspace $\goth q$ is invariant with respect to $H$ in the adjoint representation ${\rm Ad}$. It can be identified with the tangent space to $G/H$ at the point $x^0=He$.

The dimension of Cartan subspaces of $\goth  q$ (maximal Abelian subalgebras in $\goth  q$ consisting of semisimple elements) is called the rank of $G/H$.

Now let $G/H$ be a {\it symplectic} manifold. Then $\goth  h$
has a non-trivial centre $Z({\goth  h})$. For simplicity we assume that
$G/H$ is an orbit ${\rm {Ad}} \, G\cdot Z_0$ of an element $Z_0\in{\goth  g}$. In particular, then $Z_0\in Z(\goth h)$.

Further, we can also assume that $G$ is {\it simple}. Such spaces $G/H$ are divided into 4 classes (see [3],[4]):

$(a)$ \ Hermitian symmetric spaces;

$(b)$ \ semi-K\"ahlerian symmetric spaces;

$(c)$ \ para-Hermitian symmetric spaces;

$(d)$ \ complexifications of spaces of class (a).

Dimensions of $Z({\goth  h})$ are 1,1,1,2, respectively. Spaces of class $(a)$ are Riemannian, of other three classes are pseudo-Riemannian (not Riemannian).

We focus on spaces of class $(c)$. Here the center $Z(\goth h)$ is one-dimensional, so that $Z(\goth h)={\Bbb R}Z_0$, and $Z_0$ can be normalized so that the operator $I=(\rm{ad} Z_0)_{\goth  q}$ on $\goth  q$ has eigenvalues $\pm 1$. A symplectic structure on $G/H$ is defined by the bilinear form $\omega(X,Y) = B_{\goth g}(X,IY)$ on $\goth  q$.

The $\pm 1$-eigenspaces ${\goth  q}^{\pm}\subset {\goth  q}$ of $I$ are
Lagrangian, $H$-invariant, and irreducible. They are Abelian
subalgebras of $\goth  g$. So $\goth  g$ becomes a graded Lie algebra:
$$
{\goth  g} = {\goth  q}^- +{\goth  h} + {\goth  q}^+ ,
$$
with commutation relations $[\goth h, \goth h]\subset \goth h$,
$[\goth h,\goth q^-]\subset\goth q^-$, $[\goth h,\goth q^+]\subset\goth q^+$.

The pair $({ \goth q}^+,{ \goth q}^-)$ is a Jordan pair with
multiplication
$$
\{XYZ\}=\frac{1}{2} \,  [[X,Y],\,Z],
$$
see [5]. Let $r$ and $\varkappa$ be the rank and the genus of this Jordan pair. The rank $r$ coincides with the rank of $G/H$.

Set $Q^{\pm} =\exp {\goth  q}^{\pm}$. The subgroups $P^{\pm}=HQ^{\pm}=
Q^{\pm}H$ are maximal parabolic subgroups of $G$. One has the following decompositions:
\begin{eqnarray}
G&=&\overline{Q^+HQ^-} \\
&=&\overline{Q^-HQ^+},
%\eqno (1.1)(1.2)
\end{eqnarray}
where bar means closure and the sets under the bar are open and dense in
$G$. Let us call (1.1) and (1.2) the {\it{Gauss decomposition}} and (allowing some slang) the {\it{anti-Gauss decomposition}} respectively. Decompositions (1.1), (1.2) mean that almost any element $g\in G$ can be decomposed as ( the Gauss decomposition):
$$
g=\exp \eta \cdot h \cdot \exp \xi,
\eqno(1.3)
$$
or (the anti-Gauss decomposition):
$$
g=\exp \xi \cdot h \cdot \exp \eta,
\eqno(1.4)
$$
where $h\in H$, $\xi \in {\goth q^-}$, $\eta\in {\goth q^+}$, all three factors in (1.3) and (1.4) are defined uniquely. We also use the Gauss decomposition (1.3) in a little different form:
$$
g=\exp \eta \cdot \exp \xi \cdot h,
\eqno(1.5)
$$
where $\eta$ and $h$ are the same as in (1.3),  and $\xi$ is obtained from $\xi$ in (1.3) by ${\rm Ad} \, h$.

Decompositions (1.3) and (1.4) generate actions of $G$ on ${\goth  q}^{-}$ and ${\goth  q}^{+}$ respectively, namely, $\xi \mapsto\widetilde \xi=\xi\bullet g$ and $\eta \mapsto\widehat\eta=\eta\circ g$:
\begin{eqnarray}
{\rm {exp }}\, \xi \cdot g&=&{\rm {exp}}\, Y \cdot \widetilde h\cdot
{\rm {exp }}\, \widetilde \xi, \\
{\rm {exp }}\, \eta \cdot g&=&{\rm {exp}}\, X\cdot \widehat h\cdot
{\rm {exp }}\, \widehat \eta,
%\eqno (1.6)(1.7)
\end{eqnarray}
where $X\in {\goth q^-},\, Y\in {\goth q^+}$. These actions are defined on open and dense sets depending on $g$. Therefore, $G$ acts on $\goth q^-\times \goth q^+:\, (\xi ,\eta )\mapsto (\widetilde \xi ,\widehat \eta )$. The stabilizer of the point $(0,0)\in {\goth q^-\times \goth q^+}$ is $P^+\cap P^-=H$, so that we get an embedding
$$
\goth q^-\times \goth q^+\hookrightarrow G/H.
\eqno(1.8)
$$
It is defined on an open and dense set, its image is also an open and dense set. Therefore, we can consider $(\xi ,\eta)\in \goth q^-\times \goth q^+$ as coordinates on $G/H$, let us call them {\it horospherical coordinates}.

Let us write explicit formula for embedding (1.8). We use a redecomposition "anti-Gauss"\ to "Gauss". We take $\xi\in\goth q^-$, $\eta\in\goth q^+$ and decompose the anti-Gauss product ${\rm {exp}} \,\xi \cdot {\rm {exp}} \,(-\eta )$ according to formula (1.5) (the "Gauss"):
$$
{\rm {exp}} \, \xi \cdot {\rm {exp}} \, (-\eta )={\rm {exp}} \, Y \cdot  {\rm {exp}} \, X \cdot h,
\eqno(1.9)
$$
where $X\in {\goth q^-}$, $Y\in {\goth q^+}$. The obtained element $h\in H$ depends on $\xi$ and $\eta$ only, denote it by $h(\xi ,\eta )$.  Using (1.9), let us form the following element $g\in G$:
$$
g={\rm {exp}} \, Y {\rm {exp}} \, \xi =
{\rm {exp}} \, X \cdot h \cdot {\rm {exp}} \, \eta.
\eqno(1.10)
$$
Then the pair $\xi ,\eta $ goes just to the point $x=x^0 g$ where $g$ is defined by (1.10).

Under the action of the group $G$ the element $h(\xi ,\eta )$ is transformed as follows:
$$
h(\widetilde \xi ,\widehat \eta )=
{\widetilde h}^{-1}\cdot h(\xi ,\eta ) \cdot {\widehat h},
\eqno(1.11)
$$
where $\widetilde h$ and $\widehat h$ are taken from (1.6) and (1.7) respectively.

For $h\in H$, denote
$$
b(h)={\rm det \, (Ad} \, h)|_{{\goth q}^+} \, .
$$

The function $k(\xi,\eta)=b(h(\xi,\eta))$ is $N(\xi ,\eta )^{-\varkappa}$, where
$N(\xi ,\eta)$ is an irreducible polynomial $N(\xi ,\eta )$ of degree $r$ in $\xi $ and in $\eta $ separately. Considered as a function on $G/H$, the function $k(\xi,\eta)$ is an analogue of the Bergman kernel for Hermitian symmetric spaces. It follows from (1.11) that under action of $g\in G$ the function $N(\xi ,\eta)$ is transformed as follows:
$$
N(\, \widetilde \xi \, , \, \widehat \eta \,  )=
N(\xi ,\eta) \cdot b(\, {\widetilde h}\, )^{1/\varkappa}\cdot b(\, {\widehat h}\, )^{-1/\varkappa}
\eqno(1.12)
$$
In horospherical coordinates the $G$-invariant measure on $G/H$ is:
$$
dx=dx(\xi,\eta)=|N(\xi ,\eta )|^{-\varkappa} \, d\xi \, d\eta,
$$
where $d\xi$ and $d\eta$ are Euclidean measures on ${\goth q^-}$ and ${\goth q^+}$ respectively.

\null

\null
\setcounter{theorem}{0}
\setcounter{section}{2}
\setcounter{equation}{0}
\setcounter{lemma}{0}

\begin{center}
{{\large {\bf \S \ 2. Maximal degenerate series representations}}}
\end{center}

\null

In this Section we introduce two series of representations induced by characters of maximal parabolic subgroups $P^{\pm}$ of $G$ (maximal degenerate series representations).

First we take the character $\omega_\lambda (h)$, $\lambda\in \Bbb C$, of $H$:
$$
\omega_\lambda (h)=\big|b(h)\big|^{-\lambda/\varkappa}
$$
and then we extend this character to the subgroups $P^{\pm}$, setting it equal
to 1 on $Q^{\pm}$.

Let us consider induced representations $\pi^{\pm }_\lambda$ of $G$:
$$
\pi^{\pm }_\lambda ={\rm {Ind}} \, \left(G, \, P^{\mp}, \, \omega_{\mp \lambda }\right).
$$
They act on the space ${\cal {D}}^{\pm }_\lambda (G)$ of functions $f\in C^\infty (G)$ having the uniformity property:
$$
f(pg)=\omega_{\mp \lambda}(p)f(g), \ \ p\in P^{\mp },
$$
by translations from the right:
$$
\left(\pi ^{\pm }_\lambda (g)f\right)(g_1)=f(g_1g).
$$

Realize them in the {\it non-compact picture}: we restrict functions from ${\cal {D}}^{\pm }_\lambda (G)$ to the subgroups $Q^{\pm}$ and identify them (as manifolds) with ${\goth q}^{\pm}$, we obtain
\begin{equation*}
\left(\pi^-_\lambda (g)f\right)(\xi )=\omega_\lambda (\widetilde h)f
(\widetilde \xi ), \ \ \
\left(\pi^+_\lambda (g)f\right)(\eta)=\omega_{\lambda} ({\widehat h}^{-1})f
(\widehat \eta ),
\end{equation*}
where $\widetilde \xi$, $\widetilde h$, $\widehat \eta$, $\widehat h$ are taken from decompositions (1.6), (1.7).

Let us write intertwining operators. Introduce operators $A_\lambda$ and $B_\lambda$ by:
\begin{eqnarray*}
(A_\lambda \varphi )(\eta)&=&\int_{\goth q^-} |N(\xi ,\eta )|^{-\lambda-\varkappa} \
\varphi (\xi) \, d\xi, \\
(B_\lambda \varphi )(\xi)&=&\int_{\goth q^+} |N(\xi ,\eta )|^{-\lambda-\varkappa} \
\varphi (\eta) \, d\eta.
\end{eqnarray*}
The operator $A_\lambda$ intertwines $\pi^-_\lambda$ with $\pi^+_{-\lambda -\varkappa}$ and the operator $B_\lambda$ intertwines $\pi^+_\lambda$ with $\pi^-_{-\lambda -\varkappa}$ \, . Their composition is a scalar operator:
$$
B_\lambda A_{-\lambda-\varkappa } = c(\lambda)^{-1} \cdot {\rm id},
\eqno(2.1)
$$
where $c(\lambda)$ is a meromorphic function of $\lambda$, it is invariant with respect to the change $\lambda \mapsto-\lambda-\varkappa$.

We can extend $\pi^{\pm }_\lambda $, $A_\lambda$ and $B_\lambda$ to distributions on ${\goth q^-}$ and ${\goth q^+}$.

The representation $\pi ^-_\lambda$ of the Lie algebra $\goth g$ is given by some differential operators of the first order. This representation can be considered on different spaces of functions of $\xi $: for example, the space $C^\infty ({\goth q}^-)$, the space ${\rm {Pol}}\,({\goth q}^-)$ of polynomials in $\xi $, the space ${\cal {D}}'({\goth q}^-)$ of distributions on ${\goth q}^-$, in particular, the space $Z$ of distributions on ${\goth q}^-$ concentrated at the origin, etc. The same concerns to $\pi ^+_\lambda $.

Notice
\setcounter{equation}{1}
\begin{equation}
\omega_{\lambda}\big(h(\xi,\eta)\big)=|N(\xi ,\eta)|^\lambda \, ,
%\eqno(2.2)
\end{equation}
hence formula (1.12) gives
$$
|N(\widetilde \xi ,\widehat \eta )|^\lambda=
|N(\xi ,\eta)|^\lambda \cdot \omega_\lambda({\widetilde h})^{-1}
\cdot\omega_\lambda({\widehat h}) \, ,
$$
which can be interpreted as an invariance property of the function $|N(\xi ,\eta)|^\lambda $ \, :
\begin{equation}
\Big[\pi^-_\lambda (g)\otimes \pi^+_\lambda(g)\Big]
|N(\xi ,\eta)|^\lambda = |N(\xi ,\eta)|^\lambda.  \nonumber
\end{equation}

\null

\null

\begin{center}
{\large {\bf \S \ 3. Symbols and transforms }}
\end{center}

\null

In this Section we apply to a para-Hermitian symmetric space $G/H$ the scheme of quantization in the spirit of Berezin offered in [6]. We consider a variant of the quantization, we call the {\it polynomial  quantization}. For an initial algebra of operators we take here the algebra of operators $\pi^-_\lambda({\rm Env} \, (\goth g))$, where $\lambda \in \Bbb C$ and ${\rm Env} \, (\goth g)$ is the universal enveloping algebra of $\goth g$. In contrast to [6], we use the non-compact picture for representations $\pi ^{\pm }_\lambda $, see \S \ 2. The role of the Fock space is played by a space of functions $\varphi (\xi)$, $\xi \in {\goth q}^-$, so that our operators act in functions $\varphi (\xi )$. We introduce  covariant and contravariant symbols of operators, the Berezin transform etc.

As a (an analogue of) {\it supercomplete system} we take the kernel of the intertwining operators $A_{-\lambda -\varkappa}$ from \S \ 2, namely, the function
$$
\Phi(\xi ,\eta )=\Phi _{\lambda}(\xi ,\eta )=|N(\xi ,\eta)|^{\lambda} \, .
\eqno(3.1)
$$
For an operator $D=\pi^-_\lambda(X)$, $X\in{\rm Env} \, (\goth g))$, the function
$$
F(\xi,\eta)=\frac 1{\Phi (\xi,\eta)} \,
(\pi^-_\lambda (X)\otimes 1) \Phi(\xi,\eta)
$$
is called the {\it covariant symbol} of $D$. Since $\xi$, $\eta$ are horospherical coordinates on $G/H$, covariant symbols become functions on $G/H$ and, moreover, {\it polynomials} on $G/H\subset{\goth g}$. It is why we call this variant of quantization the {\it polynomial} quantization. Denote the space of covariant symbols by ${\cal {A}}_\lambda $. For generic $\lambda$, this space is the space $S(G/H)$ of {\it all polynomials} on $G/H$.

In particular, the covariant symbol of the identity operator is the function on $G/H$ equal to $1$ identically. If $X$ belongs to the Lie algebra $\goth g$ itself, then the covariant symbol of the operator $\pi^-_\lambda (X)$ is a linear function $B_{\goth g}(X,x)$ of $x\in G/H\subset {\goth g}$ with coordinates $\xi$, $\eta$, up to a factor depending on $\lambda$.

The operator $D$ is recovered by its covariant symbol $F$:
$$
(D\varphi )(\xi )=c \, \int_{G/H} \,  F(\xi ,v) \, \frac {\Phi (\xi ,v)}{\Phi (u,v)} \, \varphi (u)\,dx(u,v),
\eqno(3.2)
$$
where $c=c(\lambda)$ is taken from (2.1). Indeed, the function $\Phi$ has a reproducing property
$$
\varphi(\xi)=c(\lambda) \int_{G/H} \
\frac{\Phi(\xi,\eta)}{\Phi(u,v)} \ \varphi (u) \, dx(u,v).
$$
which is nothing but formula (2.1) written in another form.

Let $U$ be the representation of the group $G$ by translations in functions on $G/H$ (quasiregular representation), for example, on the space $C^{\infty}(G/H)$, and $U$ the corresponding representation of the Lie algebra $\goth g$. The correspondence $D\mapsto F$, assigning to an operator its covariant symbol, is $\goth g$-equivariant, it means that if $F$ is the covariant symbol of an operator $D=\pi ^-_\lambda (X),\, X\in {\rm {Env}}\,(\goth g)$, then $U(L)F$, where $L\in {\goth g}$, is the covariant symbol of the operator $\pi ^-_\lambda({\rm {ad}}\,L\cdot X)$.

The multiplication of operators gives rise to the multiplication of covariant symbols, denote it by $*$. Namely, let $F_1$, $F_2$ be covariant symbols of operators $D_1$, $D_2$ respectively. Then the covariant symbol $F_1*F_2$ of the product $D_1D_2$ is
$$
(F_1*F_2)(\xi,\eta)= \frac 1{\Phi (\xi,\eta)} \,
(D_1\otimes 1) \left(\Phi(\xi,\eta)F_2(\xi,\eta)\right) \, .
$$
Putting in (3.1) $D=D_1$, $F=F_1$ and $\varphi(u) = \Phi(u,\eta)F_2(u,\eta)$, we get
$$
(F_1*F_2)(\xi ,\eta )=\int_{G/H} \, F_1(\xi ,v)F_2(u,\eta ) \, {\cal {B}}(\xi,\eta;u,v) \, dx(u,v),
\eqno(3.3)
$$
where
$$
{\cal {B}}(\xi ,\eta ;u,v)=c \, \frac {\Phi (\xi ,v)
\Phi (u,\eta )}{\Phi(\xi ,\eta )\Phi (u,v)} \, .
$$
Let us call this function ${\cal {B}}$ the {\it Berezin kernel}. It can be regarded as a function ${\cal B}(x,y)$ on $G/H\times G/H$. It is invariant with respect to $G$:
$$
{\cal {B}}({\rm {Ad}}\,g\cdot x,\,{\rm {Ad}}\,g\cdot y)=
{\cal {B}}(x,y).
$$

\null

Define a transform $x\mapsto x \breve{\hskip 5pt} $ of the space $G/H$, that in horospherical coordinates $\xi ,\eta $ is the permutation $\xi \leftrightarrow \eta $, i.~e. $(\xi ,\eta )\mapsto (\eta ,\xi )$. This transform induces the transform $F\mapsto F\, \breve{\hskip 5pt}$ of functions in $S(G/H)$: $F\, \breve{\hskip 5pt}(\xi ,\eta )=F(\eta ,\xi )$. The Berezin kernel is invariant with respect to the simultaneous permutation $\xi \leftrightarrow \eta $ and $u\leftrightarrow v$. By (3.3) it implies that the transform $F\mapsto F \, \breve{\hskip 5pt} $ is an anti-involution with respect to the multiplication of symbols:
$$
({F_1*F_2})\breve {\hskip 5pt}= F_2\breve{\hskip 5pt} * F_1\breve{\hskip 5pt}.
$$
For operators $D$, the transform $F\mapsto F \, \breve{\hskip 5pt}$ of symbols means the conjugation $D\mapsto D \, \breve{\hskip 5pt} $ with respect to the bilinear form generated by the operator  $A_{\lambda}$:
$$
(A_{\lambda}\varphi ,\psi )=\int |N(\xi ,\eta )|^{-\lambda -\varkappa}
\varphi (\xi )\psi (\eta ) \, d\xi \, d\eta .
$$
Moreover, if $D=\pi^-_\lambda (X)$, then
$$
D \, \breve{\hskip 5pt} =\pi ^+_\lambda (X^\vee ) \, ,
$$
where $X\mapsto X^\vee$ is the transform of the algebra ${\rm {Env}}\,(\goth g)$ induced by the transform $g\mapsto g^{-1}$ in the group $G$.

Thus, the space ${\cal {A}}_\lambda$ with multiplication $*$ is an associative algebra with 1, the transform $F\mapsto F \, \breve{\hskip 5pt} $ is an anti-involution of this algebra.

\null

Now we define {\it contravariant symbols}. A function (a polynomial) $F(\xi ,\eta )$ is the contravariant symbols for the following operator $A$ (acting on functions $\varphi (\xi )$):
$$
(A\varphi )(\xi )=c \, \int_{G/H} \, F(u,v) \, \frac {\Phi (\xi ,v)}{\Phi (u,v)} \, \varphi (u) \, dx(u,v)
\eqno(3.4)
$$
(notice that (3.4) differs from (3.2) by the first argument of $F$ only). This operator is a T\"oplitz type operator.

Thus, we have two maps: $D\mapsto F$ ("co") and $F\mapsto A$ ("contra"), connecting polynomials on $G/H$ and operators acting on functions $\varphi (\xi )$.

If a polynomial $F$ on $G/H$ is the covariant symbol of an operator $D=\pi^-_\lambda (X)$, $X\in {\rm {Env}}\,(\goth g)$, and the contravariant symbol of an operator $A$ simultaneously, then
$A=\pi^-_{-\lambda-\varkappa }(X^{\vee })$. Therefore, $A$ is obtained from $D$ by the conjugation with respect to the bilinear form
$$
(F,f)=\int_{\goth q^-} \, F(\xi )\, f(\xi ) \, d\xi.
$$
In terms of kernels, it means that the kernel $L(\xi ,u)$ of the operator $A$ is obtained from the kernel $K(\xi ,u)$ of the operator $D$ by the transposition of arguments and the change $\lambda\mapsto-\lambda -\varkappa$. So, the composition ${\cal {O}}: D\mapsto A$ ("contra"\ $\circ$ "co") is
$$
{\cal {O}}: \ \pi^-_\lambda (X)\longmapsto \pi^-_{-\lambda -\varkappa }(X^{\vee }).
$$
This map commutes with the adjoint representation. Such a map was absent in Berezin's theory for Hermitian symmetric spaces.

The composition ${\cal {B}}$ ("co"\  $\circ$ "contra") maps the contravariant symbol of an operator $D$ to its covariant symbol. Let us call ${\cal {B}}$ the {\it {Berezin transform}}. The kernel of this transform is just the  Berezin kernel.

Main problems here are the following. One has to do explicitly: (a) to express the Berezin transform ${\cal {B}}$ in terms of Laplacians $\Delta_1,\ldots,\Delta_r$ ($r$ being the rank), in fact, it is the same that to decompose a canonical representation (see \S \ 4) into irreducible constituents; (b) to compute eigenvalues of ${\cal {B}}$ on irreducible constituents; (c) to find a full asymptotics of ${\cal {B}}$ when $\lambda \to -\infty$ (an analogue of the Planck constant is $h=-1/\lambda$). In particular, two first terms of asymptotic decomposition the Berezin transform ${\cal {B}}$ when $\lambda \to -\infty $ should be
$$
{\cal {B}}\sim 1-\frac 1\lambda \, \Delta ,
$$
where $\Delta $ is the Laplace--Beltrami operator. It gives the correspondence principle:
\begin{eqnarray*}
F_1*F_2 &\longrightarrow& F_1F_2,   \\
-\lambda \,(F_1*F_2-F_2*F_1)&\longrightarrow & \{F_1,\,F_2\},
\end{eqnarray*}
in right hand sides the pointwise multiplication and the Poisson bracket stand.

These problems are solved for spaces of rank one and for spaces with the group $G={\rm SO}_0(p,q)$, for latter spaces the rank is equal to 2, see [7], [8].

\null

\null

\setcounter{theorem}{0}
\setcounter{section}{4}
\setcounter{equation}{0}
\setcounter{lemma}{0}

\begin{center}
{\large\bf \S \ 4. {Canonical representations and quantization }}
\end{center}

\null

The main tool for studying quantization is the so-called {\it canonical representations} (this term was introduced in [9]). For Hermitian symmetric spaces $G/K$, they were introduced by Vershik, Gelfand and  Graev [9] (for the Lobachevsky plane) and Berezin [1], [2] (in classical case). These representations act by translations in functions on $G/K$ and are unitary with respect to some non-local inner product (now called a {\it Berezin form}).

We define canonical representations of a group $G$ in a more general setting. We give up the condition of unitarity (as too narrow) and let these representations act on sufficiently extensive spaces, in particular, on spaces of distributions. Moreover, we permit also non-transitive actions of a group $G$. Our approach uses the notion of an "overgroup"\ and consists in the following.

Let $G$ and $\widetilde G$ be semisimple Lie groups such that $G$ is a spherical subgroup of the  $\widetilde G$ (i.~e. $G$ is the fixed point subgroup of an involution of $\widetilde G$). We call $\widetilde G$ an "overgroup"\ for $G$. Let $\widetilde P$ be a maximal parabolic subgroup of $\widetilde G$, such that $\widetilde P\cap G=H$. Let ${\widetilde R}_{\lambda}$, $\lambda\in\Bbb C$, be a series of representations of $\widetilde G$ induced by characters of $\widetilde P$. They can depend on some discrete parameters, we do not write them. As a rule, representations ${\widetilde R}_{\lambda}$ are irreducible. They act on a compact manifold $\Omega$ (a flag manifold for $\widetilde G$). The series $R_{\lambda}$ has an intertwining operator $Q_{\lambda}$. Restrictions $R_{\lambda}$ of ${\widetilde R}_{\lambda}$ to $G$ are called canonical representations of $G$. They act on functions on $\Omega$. In general, $\Omega$ is not a homogeneous space for $G$, there are several open $G$-orbits on $\Omega$. They are semisimple symmetric spaces $G/H_i$. The manifold $\Omega$ is the closure of the union of open orbits. The intertwining operator $Q_{\lambda}$ called the Berezin transform.

One can consider a some different version of canonical representations, namely, the restriction of canonical representations in the first sense (on functions on $\Omega$) to some open orbit $G/H_i$. Both variants deserve to be investigated. The first variant is in some sense more natural. But for quantization we need just the second variant.

\null

Recall, see \S \ 1, classification (a), (b), (c), (d) for symplectic symmetric spaces $G/H$.  As an overgroup $\widetilde G$ for $G$ for classes (a), (b) we take the complexification $G^{ \, \Bbb C}$ of $G$, and for classes (c), (d) we take the direct product $G\times G$.

\null

\null

\setcounter{theorem}{0}
\setcounter{section}{5}
\setcounter{equation}{0}
\setcounter{lemma}{0}

\begin{center}
{\large\bf \S \ 5. {Polynomial quantization and the overgroup }}
\end{center}

\null

Now let $G/H$ be a para-Hermitian symmetric space, see \S \ 2. As an overgroup for $G$ we take the direct product ${\widetilde G}=G\times G$. It contains $G$ as the diagonal $\{(g,g), g\in G\}$.

First we describe a series of representations ${\widetilde R}_{\lambda}$ of $\widetilde G$.

Let $\widetilde P$ be a parabolic subgroup $\widetilde P$ consisting of pairs $(zh,hn)$, $z\in Q^-$, $h\in H$, $n\in Q^+$. Let ${\widetilde \omega}_\lambda$ be a character of $\widetilde P$ equal to $\omega_\lambda(h)$ at these pairs. The representation of $\widetilde G$ induced by the character ${\widetilde \omega}_{\lambda}$ of the subgroup $\widetilde P$ is denoted ${\widetilde R}_{\lambda}$.

Let us give some realizations of representations ${\widetilde R}_{\lambda}$.

Denote by ${\cal C}$ the manifold of "double"\ cosets
$$
y=s^{-1}_1 Q^-Q^+ s_2, \ \ s_1, s_2 \in G.
$$
This manifold is an analogue of a cone for representations of the pseudo-orthogonal group associated with a cone. The group $\widetilde G$ acts on ${\cal C}$ as follows:
$$
y\mapsto g^{-1}_1 y g_2, \ \ g_1, g_2\in G.
\eqno(5.1)
$$
Denote by ${\cal D}_\lambda (\cal C)$ the space of functions $f$ on ${\cal C}$ of class $C^\infty$ satisfying the following homogeneity condition
$$
f(s^{-1}_1 \, h \, Q^-Q^+ s_2)=\omega_\lambda(h) \, f(s^{-1}_1 \, Q^-Q^+ s_2).
\eqno(5.2)
$$
The representation ${\widetilde R}_{\lambda}$ acts on ${\cal D}_\lambda (\cal C)$ by
$$
\big( \, {\widetilde R}_{\lambda}(g_1, \, g_2) \,  f\big)(y)=f(g^{-1}_1 y g_2), \ \ g_1, g_2\in G.
$$

Let us take in ${\cal C}$ two sections:  "hyperbolic"\ section $\cal X$ and "parabolic"\ section $\Gamma$.

The manifold $\cal X\subset {\cal C}$ consists of cosets
$$
x=s^{-1} Q^-Q^+ s, \ \ s \in G.
$$
The group $G$ acts on $\cal X$ by $x\mapsto g^{-1}xg$. The stabilizer of the initial point $x^0=Q^-Q^+ $ is $H$, so that $\cal X$ can be identified with $G/H$.

The manifold $\Gamma\subset {\cal C}$ consists of cosets
$$
\gamma={\rm {exp}} \, (-\eta) \, Q^-Q^+  {\rm {exp}} \, \xi, \ \
\xi\in {\goth q^-},  \ \eta\in {\goth q^+}.
\eqno(5.3)
$$
This manifold can be identified with ${\goth q^-}\times{\goth q^+}$.

Let us embed $\Gamma\hookrightarrow \cal X$. It is the embedding $\goth q^-\times \goth q^+ \hookrightarrow G/H$, see (1.8), (1.9), (1.10) and further, in terms of ${\cal C}$.

Let a point $x=s^{-1} Q^-Q^+ s$, $s \in G$, has horospherical coordinates $\xi,\eta$. By (1.9) we find the element $h(\xi,\eta)$:
$$
{\rm {exp}} \, \xi \cdot {\rm {exp}} \,(-\eta )=
{\rm {exp}} \, (-Y) \cdot {\rm {exp}} \, X \cdot h_0, \ \ h_0=h(\xi,\eta),
$$
where $X\in {\goth q^-}$, $Y\in {\goth q^+}$, and by (1.10) we obtain
$$
s={\rm {exp}} \, Y \cdot {\rm {exp}} \, \xi =
{\rm {exp}} \, X \cdot h_0 \cdot {\rm {exp}} \, \eta,
\eqno(5.4)
$$
so that
\setcounter{equation}{4}
\begin{eqnarray}
x&=&s^{-1} Q^-Q^+ s \nonumber\\
&=&{\rm {exp}} \, (-\eta)\cdot h^{-1}_0 Q^-Q^+  {\rm {exp}} \, \xi.
%\eqno(5.5)
\end{eqnarray}
Thus, the embedding above assigns to a point $\gamma\in\Gamma$, given by (5.3), the point
$x\in\cal X$, given by (5.5) where $h_0=h(\xi,\eta)$.

The representation ${\widetilde R}_{\lambda}$ can be realized in functions on these manifolds  $\cal X$ and $\Gamma$.

First consider $\cal X$. A point $x=s^{-1} Q^-Q^+ s$ in $\cal X$ under action (5.1) goes to the point $g^{-1}_1 x g_2=g^{-1}_1 s^{-1} Q^-Q^+ s g_2$ in $\cal C$. Take the element $sg_2(sg_1)^{-1}$, i.~e. the element $sg_2g_1^{-1}s^{-1}$, and decompose it "by Gauss":
$$
sg_2g_1^{-1}s^{-1}={\rm exp}(-Y^*) \cdot {\rm exp} \, X^* \cdot h^*, \ X^*\in {\goth q^-}, \
Y^*\in {\goth q^+}.
\eqno(5.6)
$$
Here the element $h^*\in H$ depends on the point $x$ only and does not depend on its representative $s$. Let us form an element $s^*\in G$:
$$
 s^*={\rm exp} \, Y^* \cdot sg_2 ={\rm exp} \, X^* \cdot h^*sg_1.
\eqno(5.7)
$$
It gives the point $x^*\in \cal X$:
$$
 x^*=(s^*)^{-1} Q^-Q^+ s^*.
$$
By (5.7) we have
$$
x^*=g_1^{-1}s^{-1} (h^*)^{-1}Q^-Q^+ sg_2.
$$
Therefore,
$$
f( x^*)=\omega_\lambda((h^*)^{-1}) \, f(g_1^{-1}x g_2),
$$
so that $ {\widetilde R}_{\lambda}$ acts in functions on ${\cal X}=G/H$ as follows:
$$
\big({\widetilde R}_{\lambda} \, (g_1,g_2) \, f\big)(x)= \omega_\lambda(h^*) \, f(x^*).
\eqno(5.8)
$$

\begin{theorem}
In horospherical coordinates $\xi,\eta$ on $G/H$ the representation ${\widetilde R}_{\lambda}$ is
$$
\big({\widetilde R}_{\lambda} \, (g_1,g_2) \, f \big)(\xi, \eta)=
\frac{\Phi_\lambda(\xi\bullet g_2, \eta\circ g_1)}
{\Phi_\lambda(\xi, \eta)} \
  \omega_\lambda({\widetilde h}_2) \ \omega_\lambda({\widehat h}_1^{-1}) \
f(\xi\bullet g_2, \, \eta\circ g_1),
\eqno(5.9)
$$
where $\widetilde h_2$ and $\widehat h_1$ are taken from decompositions $(1.6)$ and $(1.7)$ with $g=g_2$ and $g=g_1$ respectively.
\end{theorem}
{\bf Proof.} Let a point $x=s^{-1} Q^-Q^+ s$, $s \in G$, has horospherical coordinates $\xi,\eta$. By (5.4) and (1.6), (1.7) we have
\setcounter{equation}{9}
\begin{eqnarray}
&&sg_2={\rm {exp}} \, Y \cdot {\rm {exp}} \, \xi \,\cdot g_2=
{\rm {exp}} \, Y \cdot {\rm {exp}} \, Y_2 \,\cdot {\widetilde h}_2 \cdot
{\rm {exp}} \, {\widetilde \xi}_2, \nonumber\\
&&sg_1={\rm {exp}} \, X \cdot h_0 \cdot {\rm {exp}} \,\eta \, \cdot g_1=
{\rm {exp}} \, X \cdot h_0 \cdot{\rm exp} \, X_1 \cdot {\widehat h}_1 \cdot
{\rm {exp}} \, {\widehat \eta}_1. \nonumber
\end{eqnarray}
where ${\widetilde \xi}_2=\xi\bullet g_2$, ${\widehat \eta}_1=\eta\circ g_1$. Hence
\begin{eqnarray}
&&s^*={\rm {exp}} \, Y^*\cdot sg_2=
{\rm {exp}} \, Y_3 \, \cdot {\widetilde h}_2 \cdot
{\rm {exp}} \, {\widetilde \xi_2}, \\
&&s^*={\rm {exp}} \, X^* \cdot sg_1= {\rm {exp}} \, X_3
\cdot h^* \cdot h_0 \cdot {\widehat h}_1 \cdot
{\rm {exp}} \, {\widehat \eta}_1.
%\eqno(5.10),\eqno(5.11)
\end{eqnarray}
Therefore, using (5.10) and (5.11), we obtain
\begin{eqnarray*}
x^*&=&(s^*)^{-1} Q^-Q^+ s^*   \\
&=& {\rm {exp}} \, {\widehat \eta}_1 \cdot \left(h^* h_0 {\widehat h}_1\right)^{-1} \cdot
     Q^-Q^+ \cdot {\widetilde h}_2 \cdot {\rm {exp}} \, {\widetilde \xi_2}   \\
&=& {\rm {exp}} \, {\widehat \eta}_1 \cdot \left(h^* h_0 {\widehat h}_1\right)^{-1} \cdot
     {\widetilde h}_2 \cdot Q^-Q^+ \cdot{\rm {exp}} \, {\widetilde \xi_2}.
\end{eqnarray*}
By homogeneity condition (5.2) we have
$$
f(x^*)=f\left({\rm {exp}} \, {\widehat \eta}_1 \cdot Q^-Q^+ \cdot{\rm {exp}} \, {\widetilde \xi_2}\right)
\cdot  \omega_\lambda \left((h^* h_0 {\widehat h}_1)^{-1}{\widetilde h}_2 \right)
\eqno(5.12)
$$
On the other hand, by (5.5) we can write the point $x^*$ in the following form:
$$
x^*={\rm {exp}} \, {\widehat \eta}_1 \cdot (h^*_0)^{-1} Q^-Q^+ \cdot{\rm {exp}} \, {\widetilde \xi_2},
$$
where $h^*_0=h\left({\widetilde \xi}_2, {\widehat \eta}_1\right)$. Whence again by homogeneity condition (5.2) we obtain
$$
f(x^*)=f \left({\rm {exp}} \, {\widehat \eta}_1 \cdot Q^-Q^+ \cdot{\rm {exp}} \, {\widetilde \xi_2}\right)
\cdot  \omega_\lambda \left((h^*_0)^{-1}\right).
\eqno(5.13)
$$
Comparing (5.12) and (5.13) we get
$$
\omega_\lambda \left( {\widehat h}_1^{-1} h^{-1}_0 (h^*)^{-1} {\widetilde h}_2 \right)=
\omega_\lambda \left((h^*_0)^{-1}\right),
$$
whence
$$
\omega_\lambda(h^*) =\frac{\omega_\lambda(h^*_0)}{\omega_\lambda(h_0)} \,
\omega_\lambda ( {\widehat h}_1^{-1}) \, \omega_\lambda ({\widetilde h}_2).
$$
Substitute it to (5.8) and remember (2.2) and (3.1), as result we obtain (5.9). \hfill  $\square$

\null

Similarly, the representation ${\widetilde R}_{\lambda}$ can be realized in functions on the manifold
$\Gamma$, it is given by:
$$
\left({\widetilde R}_{\lambda}(g_1,g_2) \, f \right) \, (\xi, \eta)=
  \omega_\lambda({\widetilde h}_2) \ \omega_\lambda({\widehat h}_1^{-1}) \
f(\xi\bullet g_2, \eta\circ g_1).
$$
It shows that ${\widetilde R}_{\lambda}$ is equivalent to a tensor product:
$$
{\widetilde R}_{\lambda}(g_1,g_2) = \pi^-_\lambda (g_2) \otimes \pi^+_\lambda(g_1).
$$

\null

The group $\widetilde G$ contains three subgroups isomorphic to $G$. The first one is the diagonal consisting of pairs $(g,g)$, $g\in G$. The restriction of the representation ${\widetilde R}_{\lambda}$ to this subgroup is the representation $U$ by translations on $G/H$:
$$
\left({\widetilde R}_{\lambda}(g,g) \, f\right)(x)= f(g^{-1}xg).
$$
Indeed, (5.6) and (5.7) with $g_1=g_2=g$ give $h^*=e$ and $s^*=sg$.

\null

Two other subgroups $G_1$ and $G_2$ consist of pairs $(g,e)$ and $(e,g)$, where $g\in G$, respectively.

By virtue of Theorem 5.1, the restriction of the representation ${\widetilde R}_{\lambda}$ to the subgroup $G_2$ is given by
\begin{eqnarray*}
\left({\widetilde R}_{\lambda} (e,g) \, f \right)(\xi, \eta)&=&
\frac{\Phi_\lambda({\widetilde\xi}, \eta)}
{\Phi_\lambda(\xi, \eta)} \
  \omega_\lambda({\widetilde h}) \, f({\widetilde\xi}, \eta)  \\
&=&\frac{1}{\Phi_\lambda(\xi, \eta)} \left(\pi^-_\lambda (g)\otimes 1\right)
\Big[ f(\xi, \eta) \Phi_\lambda(\xi, \eta)\Big] \, .
\end{eqnarray*}
Similarly, the restriction of the representation $ {\widetilde R}_{\lambda} $ to the subgroup $G_1$ is given by
$$
\big({\widetilde R}_{\lambda} (g,e) \, f \big)(\xi, \eta)=
\frac{1}{\Phi_\lambda(\xi, \eta)} \, \big(1\otimes \pi^+_\lambda(g)\big)
\Big[ f(\xi, \eta) \Phi_\lambda(\xi, \eta)\Big] \, .
$$

Let us go from the group $G$ to the universal enveloping algebra ${\rm Env} (\goth g)$ and preserve notations for representations. Let us take as $f$ the function $f_0$ equal to the 1 identically. Then for $X\in {\rm Env} (\goth g)$ we have
\setcounter{equation}{13}
\begin{eqnarray}
({\widetilde R}_{\lambda} (0,X)f_0)(\xi,\eta)&=&\frac{1}{\Phi_\lambda (\xi,\eta)}
(\pi^-_\lambda(X)\otimes 1)\Phi_\lambda (\xi,\eta), \\
({\widetilde R}_{-\lambda-\varkappa}(X,0)f_0)(\xi,\eta)&=&\frac{1}{\Phi_\lambda (\xi,\eta)}
(1\otimes \pi^+ _{-\lambda-\varkappa}(X))\Phi_{-\lambda-\varkappa} (\xi,\eta).
%\eqno(5.14)\eqno(5.15)
\end{eqnarray}
Right hand sides of formulae (5.14) and (5.15) are just covariant and contravariant symbols of operator  $D=\pi^-_\lambda (X)$ in polynomial quantization, see \S \ 3.

\null

Let us change the position of arguments in ${\widetilde R}_{\lambda}$, then we have a new representation ${\widehat R}_\lambda$ of ${\widetilde G}$, namely, ${\widehat R}_\lambda (g_1, g_2) = {\widetilde R}_{\lambda} (g_2, g_1)$. Using the realization of ${\widetilde R}_{\lambda}$ on the section $\Gamma$, we see that the tensor product $A_{\lambda}\otimes B_{\lambda}$ intertwines the representation ${\widetilde R}_{\lambda}$ with the representation ${\widehat R}_{-\lambda-\varkappa}$. Passing from $\Gamma$ to $\cal X$ and replacing $\lambda$ by $-\lambda-\varkappa$, we obtain that the operator $c(\lambda) A_{-\lambda-\varkappa}\otimes B_{-\lambda-\varkappa}$ intertwines the representation ${\widehat R}_{-\lambda-\varkappa}$ with the representation ${\widetilde R}_{\lambda}$ and transfers contravariant symbols to covariant ones. It has the kernel ${\cal B}_\lambda(\xi,\eta;u,v)$, i.~e. it is precisely the Berezin transform.

\null

\null

\setcounter{theorem}{0}
\setcounter{section}{6}
\setcounter{equation}{0}
\setcounter{lemma}{0}

\begin{center}
{\large\bf \S \ 6. {Polynomial quantization on rank one spaces }}
\end{center}

\null
 
We consider here the spaces $G/H$, where $G={\rm {SL}}\,(n,\mathbb R)$, $H={\rm {GL}}\,(n-1,\mathbb R)$. They have dimension $2n-2$, rank $r=1$ and genus $\varkappa=n$. These spaces $G/H$ exhaust all para-Hermitian symmetric spaces of rank one up to the covering. Further we assume $n\geqslant 3$. 

Let ${\rm Mat} \, (n,\mathbb R)$ denote the space of real $n\times n$ matrices $x$. The Lie algebra $\mathfrak g$ of $G$ consists of $x$ with ${\rm tr} \, x =0$. By \S \ 1, the space $G/H$ is a $G$-orbit in $\mathfrak g$.   

But now it is more convenient for us to change a little the realization of $G/H$. 

The group $G$ acts on ${\rm {Mat}}\,(n,\mathbb R)$ by $x\mapsto g^{-1}xg$. Let us write matrices $x$ in the block form according to the partition $n=(n-1)+1$:
$$
x=\left(
\begin{array}{cc}
\alpha&\beta \\
\gamma&\delta
\end{array}
\right)
$$
where $\alpha\in{\rm {Mat}}\,(n-1,\mathbb R)$, $\beta$ is a vector-column in ${\mathbb R}^{n-1}$, $\gamma$ is a vector-row in ${\mathbb R}^{n-1}$ and $\delta$ is a number. 
  
Let $x^0$ be the following matrix: 
$$
x^0=\left(
\begin{array}{cc}
0&0\\
0&1
\end{array}
\right)
$$
The $G$-orbit of $x^0$ is just $G/H$. This manifold is the set of matrices $x$ whose trace and rank are equal to 1. The stabilizer $H$ of $x^0$ consists of matrices ${\rm diag}\{a,b\}$, where $a\in{\rm {GL}}\,(n-1,\mathbb R)$, $b=({\det} \ a)^{-1}$, so that $H={\rm {GL}}\,(n-1,\mathbb R)$. 

Subalgebras ${\mathfrak q}^-$ and ${\mathfrak q}^+$ consist respectively of matrices
$$
X=\left(
\begin{array}{cc}
0 & 0 \\
\xi & 0
\end{array}
\right), \ \ \ 
Y=\left(
\begin{array}{cc}
0 & \eta \\
0 & 0
\end{array}
\right), 
$$
where $\xi$ is a row $(\xi_1,\ldots,\xi_{n-1})$, and $\eta$ is a column $(\eta_1,\ldots,\eta_{n-1})$ in ${\mathbb R}^{n-1}$. Embedding (1.5) is 
$$
x=\frac {1}{N(\xi,\eta)} \, \left(
\begin{array}{cc}
-\eta \xi &-\eta \\
\xi &1
\end{array}
\right),
$$
where $N(\xi,\eta)=1-\xi\eta=1-(\xi_1\eta_1+\ldots+\xi_{n-1}\eta_{n-1})$.

A $G$-invariant metric $ds^2$ on $G/H$ up to a factor is ${\rm {tr}}\,(dx^2)$. It generates the measure $dx$, the Laplace--Beltrami operator $\Delta$ and the Poisson bracket $\{f,h\}$. In coordinates $\xi,\eta$ we have:
\begin{eqnarray}
ds^2&=&-2N(\xi ,\eta )^{-2}\Bigl\{\sum\xi_i\,d\eta_i 
\sum \eta_i d\xi_i + N(\xi,\eta) \sum d\xi_i \, d\eta_i \Bigr\}. \nonumber \\
dx&=&\left|N(\xi ,\eta )\right|^{-n}\,d\xi \, d\eta \ \ \ (d\xi =
d\xi _1...d\xi _{n-1}), \nonumber \\
\Delta&=&N(\xi ,\eta ) \, \sum \, (\delta _{ij}-\xi _i\eta _j) \,\frac{\partial^2}{\partial \xi_i\,\partial \eta_j} \, , \nonumber \\
\{f,h\}&=&N(\xi,\eta ) \, \sum \, (\delta _{ij}-\xi_i\eta_j)
\left(\frac {\partial f}{\partial \eta _i}\,\frac {\partial h}
{\partial \xi_j}-\frac {\partial f}{\partial \xi_i} \, 
\frac{\partial h}{\partial \eta_j} \right). \nonumber 
\end{eqnarray}

The Berezin kernel is  
$$
{\mathcal B}(x,y)= c(\lambda) \, \frac {\Phi (\xi ,v)
\Phi (u,\eta )}{\Phi(\xi ,\eta )\Phi (u,v)} \, = \, c(\lambda) \, |{\rm tr} (xy)|^{\lambda},
$$
where
$$
c(\lambda) =\left\{2^{n+1}\pi^{n-2}\Gamma(-\lambda-n+1)
\Gamma(\lambda+1)\left[\cos\left(\lambda+\frac{n}{2}\right)\pi-\cos\frac{n\pi}{2}\right]\right\}^{-1}.
$$ 

The Berezin transform is written in terms of the Laplace--Beltrami operator $\Delta$ as follows  
$$
{\mathcal B}=\frac{\Gamma(-\lambda+\sigma) \, \Gamma(-\lambda-\sigma-n+1)}
{\Gamma(-\lambda) \, \Gamma(-\lambda-n+1)} \, ,
\eqno (6.1)
$$
the right-hand side should be regarded as a function of $\Delta=
\sigma(\sigma+n-1)$. 

\null

Now let $\lambda \to -\infty $. Then (6.1) gives
$$
{\mathcal {B}} \sim 1-\frac 1{\lambda} \, \Delta. 
$$

Hence we have
$$
F_1*F_2 \sim F_1F_2-\frac {1}{\lambda} \, N^2 \, \frac {\partial F_1}{\partial \xi } \, \frac {\partial F_2}{\partial \eta },
$$
so that for $\lambda \to -\infty$ we have
$$
F_1*F_2\longrightarrow F_1F_2, 
\eqno(6.2)
$$
$$
-\lambda \,(F_1*F_2-F_2*F_1)\longrightarrow \{F_1,\,F_2\},
\eqno(6.3)
$$
in the right hand sides of (6.2) and (6.3) the pointwise multiplication and the Poisson bracket stand, respectively. Relations (6.2) and (6.3) show that for the family of algebras of covariant symbols the {\it{correspondence principle}} is true. As the Planck constant, one has to take $h=-1/\lambda$.

Moreover, we can write not only two terms of the asymptotics but also a full asymptotic decomposition (a deformation decomposition) of ${\mathcal {B}}$ explicitly. In order to have a transparent formula, one has to expand not in powers of $h=-1/\lambda$ but use "generalized powers" of $-\lambda-n$. Then decomposition turns out to be a series terminating on any irreducible subspace of polynomials on $G/H$. 

Namely, we have the following decomposition of the Berezin transform: 
$$
{\mathcal {B}}=\sum_{k=0}^{\infty } \frac{1}{k!} \cdot \frac {\Delta \,\left[\Delta -1\cdot
n\right] \, \left[\Delta -2\cdot (n+1)\right] \, \ldots \, \left[\Delta -(k-1)(k-2+n)\right]} {(-\lambda-n)^{(k)}} \, , 
$$
where
$$
a^{(m)}=a(a-1)\ldots(a-m+1).
$$

\null

\null

\setcounter{theorem}{0}
\setcounter{section}{7}
\setcounter{equation}{0}
\setcounter{lemma}{0}

\begin{center}
{\large\bf \S \ 7. {Group quantization }}
\end{center}

\null

The group $G={\rm SL}(2,\mathbb R)$ consists of real matrices 
$$
g=\left(
\begin{array}{cc}
\alpha&\beta \\
\gamma&\delta \\
\end{array}
\right), \ \ \ \alpha\delta-\beta\gamma=1.
$$

We denote:
$$
t^{\lambda, \varepsilon}=|t|^\lambda ({\rm sgn} t)^\varepsilon, \
\  \lambda\in \mathbb C, \ \ \varepsilon\in\{0,1\}.
$$
In this section we use a little wider set of parameters: instead of $\lambda$ we take the pair $(\sigma, \varepsilon)$ where $\sigma \in \mathbb C$, $\varepsilon =0,1$,
Let us denote by ${\cal {D}}_{\sigma ,\varepsilon }(\mathbb R)$ the space of functions $f\in C^\infty (\mathbb R)$ such that the function 
$$
\widehat f(t)=t^{2\sigma,\varepsilon} f\left(1/t\right)
$$
belongs to $C^\infty (\mathbb R)$ too. The representation  $T_{\sigma, \varepsilon }$ of the group $G$ acts on ${\cal {D}}_{\sigma,\varepsilon }(\mathbb R)$ by:
\begin{equation}\label{1}
\left(T_{\sigma ,\varepsilon }(g)f\right)(t)=f\left(\frac{\alpha
t+ \gamma}{\beta t + \delta} \right) (\beta t+\delta)^{2\sigma,\varepsilon}.
\end{equation}
Let us denote by $\widehat T_{\sigma,\varepsilon}$ the contragredient representation $g\mapsto T_{\sigma,\varepsilon
}(\widehat g)$, where  
$$
\widehat g=\left(
\begin{array}{cc}
\delta&\gamma \\
\beta&\alpha \\
\end{array}
\right),
$$
so that
\begin{equation}\label{2}
\left(\widehat T_{\sigma ,\varepsilon
}(g)f\right)(t)=f\left(\frac{\delta t + \beta}{\gamma t +
\alpha}\right) (\gamma t+\alpha )^{2\sigma,\varepsilon}.
\end{equation}
Representations $T_{\sigma,\varepsilon}$ and $\widehat
T_{\sigma,\varepsilon}$ are equivalent by means of the operator 
$f\mapsto\widehat f$.

\null

An operator $A_{\sigma,\varepsilon}$, defined by formula:
$$
(A_{\sigma ,\varepsilon }f)(t)=\int_{-\infty}^{\infty }
N^{-2\sigma-2,\,\varepsilon }f(s) \, ds, \ \ \  N=N(t,s) =
1-\xi\eta,
$$
intertwines $T_{\sigma,\varepsilon}$ and ${\widehat
T}_{-\sigma-1,\varepsilon}$:
$${\widehat T}_{-\sigma-1,\varepsilon}(g) A_{\sigma,\varepsilon}=
A_{\sigma,\varepsilon}T_{\sigma,\varepsilon}(g) \, ,
$$
and also ${\widehat T}_{\sigma,\varepsilon}$ and
$T_{-\sigma-1,\varepsilon}$.

\null

Consider in $\mathbb R^3$ the bilinear form
$$
[x,y]=-x_1y_1+x_2y_2+x_3y_3.
$$
Let ${\cal {X}}$ denotes the hyperboloid of one sheet $[x,x]=1$. Realise it as the set of matrices
$$
x=\frac{1}{2}\left(
\begin{array}{cc}
1-x_3 & x_2-x_1\\
x_2+x_1 & 1+x_3 \\
\end{array}
\right)
$$
with $\det x=0$. The group $G$ acts on them by conjugations $x\mapsto g^{-1}xg$ transitively. 

Let us introduce on ${\cal {X}}$ {\it horospherical coordinates} $\xi, \eta$:
$$
x=\frac{1}{N} \, (\xi +\eta, \, \xi -\eta, \, 1+\xi\eta), \ \
$$
so that:
\begin{equation}\label{3}
x=\frac {1}{N} \, \left(
\begin{array}{cc}
-\eta \xi &-\eta \\
\xi &1
\end{array}
\right).
\end{equation}
Let us define for variables $\xi, \eta$ two vectors -- vector-line and vector-column:
$$
u=u(\xi)=\left(
\begin{array}{cc}
 \xi & 1
\end{array}
\right), \ \ \ v=v(\eta)= \left(
\begin{array}{c}
-\eta  \\
  1
\end{array}
\right)
$$
Then
\begin{equation}\label{4}
x=\frac{vu}{uv}, \ \ \ N=uv.
\end{equation}

\null

As an initial class of operators we take operators  $T_{\sigma,\varepsilon}(g)$, $g\in G$. As a supercomplete system we take the kernel of the intertwining operator $A_{-\sigma-1,\, \varepsilon}$, namely, the function 
\begin{equation}\label{5}
\Phi (\xi ,\eta )=\Phi _{\sigma ,\varepsilon }(\xi ,\eta )=
N(\xi,\eta)^{2\sigma, \, \varepsilon}
\end{equation}
of two variables $\xi,\eta$.

For an operator $T_{\sigma,\varepsilon}(g)$, its {\it {covariant symbol}} and {\it {contravariant symbol}} are the following functions respectively
\begin{align}\label{6}
F_g(\xi,\eta )&=\frac {1}{\Phi _{\sigma ,\varepsilon }(\xi ,\eta )}
\, (T_{\sigma, \varepsilon} (g)\otimes 1) \, \Phi _{\sigma
,\varepsilon }(\xi ,\eta ), \\
\label{7}
F^{\natural}_g(\xi,\eta )&=\frac {1}{\Phi _{-\sigma-1 ,\varepsilon
}(\xi ,\eta )} \, (1\otimes T_{-\sigma-1, \varepsilon} (\widehat g
\, )) \, \Phi _{-\sigma-1 ,\varepsilon }(\xi ,\eta ),
\end{align}
Consider $\xi ,\eta $ as horospherical coordinates on ${\cal
{X}}$. Them functions \eqref{6}, \eqref{7} become functions on ${\cal {X}}$.

\begin{theorem}
The covariant and contravariant symbols of the operator $T_{\sigma,\varepsilon}(g)$, $g\in~G$, are the following functions on the hyperboloid $\cal{X}$ respectively:
\begin{align}
F_g(x)&=({\rm {tr}}(xg))^{2\sigma, \varepsilon} \label{8} \\
&= \left(\frac{ugv}{uv}\right)^{2\sigma, \varepsilon} ,\label{9}\\
F^{\natural}_g(x)&={\rm {tr}}(g^{-1}x)^{-2\sigma-2, \varepsilon}  \label{10}\\
&= \left(\frac{ug^{-1}v}{uv}\right)^{-2\sigma-2, \varepsilon}. \label{11}
\end{align}
\end{theorem}

{\bf Proof.} First we prove \eqref{8}, \eqref{9}. By \eqref{6}, \eqref{5}, \eqref{1} we obtain:
\begin{align}
F_g(\xi,\eta )&=\frac {1}{N^{2\sigma ,\varepsilon }} \, \left( 1-
\frac{\alpha \xi + \gamma}{\beta \xi + \delta} \, \eta
\right)^{2\sigma, \varepsilon}
\left(\beta \xi + \delta\right)^{2\sigma, \varepsilon}  \nonumber \\
&=\left( \frac{\beta \xi + \delta - \alpha \xi \eta -
\gamma\eta}{N} \right)^{2\sigma, \varepsilon} \label{12}
\end{align}
It is easy to check, the nominator in \eqref{12} is precisely 
$u(\xi)gv(\eta)$, it together with (7.4) proves (7.9). Formula \eqref{8} follows from \eqref{12} and \eqref{3}.

Now let us prove \eqref{10}, \eqref{11}. By \eqref{7}, \eqref{5}, \eqref{2} we obtain:
\begin{align}
F^\natural_g(\xi,\eta )&=\frac {1}{N^{-2\sigma-2 ,\varepsilon }}
\, \left( 1- \xi \, \frac{\delta\eta +\beta}{\gamma \eta +  \alpha
} \right)^{-2\sigma-2, \varepsilon}
\left(\gamma\eta + \alpha\right)^{-2\sigma-2, \varepsilon}  \nonumber \\
&=\left( \frac{\gamma\eta + \alpha- \delta\xi\eta - \beta\xi}{N}
\right)^{-2\sigma-2, \varepsilon} \label{13}
\end{align}
Since
$$
g^{-1}=\left(
\begin{array}{cc}
\delta &-\beta \\
-\gamma&\alpha \\
\end{array}
\right),
$$
then the nominator \eqref{12} is precisely $u(\xi)g^{-1}v(\eta)$, it together with \eqref{4} proves \eqref{11}. Formula \eqref{10} follows from \eqref{13} and \eqref{3}. \hfill  $\square$ 

\null

\null

\begin{center}
{\large\bf  {References }}
\end{center}

\null

   1. F.~A.~Berezin. Quantization on complex symmetric spaces. Izv. Akad. Nauk SSSR, Ser. mat., 1975, vol. 39, No. 2, 363--402. English transl.: Math. USSR-Izv., 1975, vol. 9, 341--379.
\vskip 5pt
   2. F.~A.~Berezin. A connection between the co- and the contravariant symbols of operators on classical complex symmetric spaces. Dokl. Akad. Nauk SSSR, 1978, vol. 19, No. 1, 15--17. English transl.: Soviet Math. Dokl., 1978, vol. 19, No. 4, 786--789.
\vskip 5pt
   3. S.~Kaneyuki. On orbit structure of compactifications of parahermitian sym\-metric spaces. Japan. J. Math., 1987, vol. 13, No. 2, 333--370.
\vskip 5pt
   4. S.~Kaneyuki, M.~Kozai. Paracomplex structures and affine symmetric spaces. Tokyo J. Math., 1985, vol. 8,
No. 1, 81--98.
\vskip 5pt
   5. O.~Loos. Jordan Pairs. Lect. Notes in Math., 1975, {\it 460}.
\vskip 5pt
   6. V.~F.~Molchanov. Quantization on para-Hermitian symmetric spaces, Amer. Math. Soc. Transl, Ser. 2, 1996, vol. 175 (Adv. in the Math. Sci.--31), 81--95.
\vskip 5pt
   7. V.~F.~Molchanov, N.~B.~Volotova. Polynomial quantization on rank one para-Hermitian symmetric spaces. Acta Appl. Math., 2004, vol. 81, Nos. 1--3, 215--232.
   \vskip 5pt

\end{document}